\documentclass[12pt]{article}
 
\usepackage{babel}
\usepackage{isolatin1,times,a4wide,doublespace,amsmath,latexsym,amssymb}
\newtheorem{theorem}{Theorem}[section]
\newtheorem{lemma}[theorem]{Lemma}
\newtheorem{proposition}[theorem]{Proposition}

\newcommand{\Tr}{\mathrm{Tr}}
\newcommand{\Span}{\mathrm{span}}

\newcommand{\GL}{\mathrm{GL}}
\newcommand{\GAP}{\mathrm{GAP}}
\newcommand{\gap}{\mathrm{gap}}
\newtheorem{definition}[theorem]{Definition}
\newtheorem{conjecture}[theorem]{Conjecture}
\newtheorem{corollary}[theorem]{Corollary}
\newtheorem{example}[theorem]{Example}

\newcommand{\bproof}{\noindent{\bf Proof: }}
\newcommand{\eproof}{\hfill $\Box$\\}
\newcommand{\bremark}{\noindent{\bf Remark: }}
\newcommand{\eremark}{\hfill \\}
\newcommand{\w}{\omega}
\newcommand{\cM}{{\mathcal M}}
\newcommand{\e}{\varepsilon}
\newcommand{\cU}{{\mathcal U}}

\newcommand{\cA}{{\mathcal A}}
\newcommand{\cF}{{\mathcal F}}

\newcommand{\cS}{{\mathcal S}}
\newcommand{\cB}{{\mathcal B}}
\newcommand{\bN}{{\mathbb N}}
\newcommand{\bR}{{\mathbb R}}

\begin{document}
\title{Embeddings of Banach Spaces into Banach Lattices and the
Gordon-Lewis Property}

\author{P.G.\ Casazza\footnote{Supported by SNF grant DMS 970618} 
\and N.J.\ Nielsen\footnote{Supported in part by the Danish
Natural Science Research Council, grants 9503296 and 9600673}}
\date{}
\maketitle

\begin{abstract}
In this paper we first show that if $X$ is a Banach space and $\alpha$ 
is a left invariant crossnorm on $\ell_\infty\otimes X$, then there 
is a Banach lattice $L$ and an isometric embedding $J$ of $X$ into $L$,
so that $I\otimes J$ becomes an isometry of $\ell_\infty\otimes_\alpha X$
onto $\ell_\infty\otimes_m J(X)$. Here $I$ denotes the identity operator 
on $\ell_\infty$ and $\ell_\infty\otimes_m J(X)$ the canonical lattice 
tensor product. This result is originally due to G.\ Pisier (unpublished), 
but our proof is different. We then use this to characterize the Gordon-Lewis
property $\GL$ in terms of embeddings into Banach lattices. Also other structures
related to the $\GL$ are investigated. 
\end{abstract}

\section*{Introduction}

In this paper we investigate embeddings of Banach spaces into Banach lattices, 
which preserve a certain tensorial structure given a priori. This is then used 
to characterize the Gordon-Lewis property $\GL$ and related structures in Banach 
spaces.

Our basic result states that if $X$ is a Banach space and $\alpha$ is a left 
tensorial crossnorm on $\ell_\infty\otimes X$ (see Section \ref{sec0} for 
the definition), then there exist a Banach lattice $L$ and an isometric embedding
$J$ of $X$ into $L$ so that $I\otimes J$ becomes an isometric embedding of 
$\ell_\infty\otimes_\alpha X$ onto $\ell_\infty\otimes_m J(X)$. Here $I$ denotes the 
identity operator on $\ell_\infty$ and $\ell_\infty \otimes_m X$ the canonical lattice
tensor product. This result was originally proved by Pisier \cite{P1} (unpublished), 
but our construction of the Banach lattice $L$ is quite different from his. It is a 
modification of a construction given by the second named author and presented at a 
conference in Columbia, Missouri in $1994$ and is based on our Theorem \ref{thm1.5} 
below.

This result is then used to prove that a Banach space $X$ has $\GL_2$ if and only
if it embeds into a Banach lattice $L$, so that every absolutely summing operator 
from $X$ to a Hilbert space extends to an absolutely summing operator defined on
$L$. In a similar manner we prove that $X$ has the general $\GL$-propery if and
only if it embeds into a Banach lattice $L$, so that every absolutely summing
operator from $X$ to an arbitrary Banach space $Y$ extends to a cone-summing operator
from $L$ to $Y$. Some related structures in Banach spaces, e.g. the Gaussian average
property defined in \cite{CN}, are also characterized in terms of embeddings into 
Banach lattices.

In Section \ref{sec1} of the paper we investigate left tensorial crossnorms 
and prove the main result mentioned above. Section \ref{sec2} is devoted to 
the characterizations of the $\GL$-property, while Section \ref{sec3} contains 
some further applications to $\GL$-subspaces of Banach lattices of finite cotype.

Let us finally mention that L.B.\ McClaran \cite{MC} has used Pisier's result to 
characterize subspaces of quotients of Banach lattices.

\setcounter{section}{-1}

\section{Notation and Preliminaries}
\label{sec0}

In this paper we shall use the notation and terminology commonly used in
Banach space theory as it appears in \cite{LT1}, \cite{LT2} and
\cite{T}. $B_X$ shall always denote the closed unit ball of the Banach
space $X$.

If $X$ and $Y$ are Banach spaces, $B(X,Y)$ ($B(X)=B(X,X)$) denotes the
space of bounded linear operators from $X$ to $Y$ and throughout the
paper we shall identify $X\otimes Y$ with the space of $\w^*$-continuous
finite rank operators from $X^*$ to $Y$ in the canonical
manner. Further, if $1\le p<\infty$ we let $\Pi_p(X,Y)$ denote the space
of $p$-summing operators from $X$ to $Y$ equipped with the $p$-summing
norm $\pi_p$; $I_p(X,Y)$ denotes the space of all $p$-integral operators
from $X$ to $Y$ equipped  with the $p$-integral norm $i_p$ and
$N_p(X,Y)$ denotes the space of all $p$-nuclear operators from $X$ to
$Y$ equipped with the $p$-nuclear norm $\nu_p$. $X\otimes_\pi Y$ denotes
the completion of $X\otimes Y$ under the largest tensor norm $\pi$ on
$X\otimes Y$.

We recall that if $1\le p\le\infty$ then an operator $T\in B(X,Y)$ is
said to factor through $L_p$ if it admits a factorization $T=BA$,
where $A\in B(X,L_p(\mu))$ and $B\in B(L_p(\mu),Y)$ for some measure
$\mu$ and we denote the space of all operators which factor through
$L_p$ by $\Gamma_p(X,Y)$. If $T\in\Gamma_p(X,Y)$ then we define
\[
\gamma_p(T) = \inf\{\|A\|\, \|B\|\, \mid \, T=BA,\, \mbox{$A$ and $B$
as above$\}$};
\]
$\gamma_p$ is a norm on $\Gamma_p(X,Y)$ turning it into a Banach
space. All these spaces are operator ideals and we refer to the above
mentioned books and \cite{HK}, \cite{KW} and \cite{Pi} for further
details. To avoid misunderstanding we stress that in this paper a
$p$-integral operator $T$ from $X$ to $Y$ has a $p$-integral
factorization ending in $Y$ with $i_p(T)$ defined accordingly; in some
books this is referred to as a strictly $p$-integral operator.

If $(\cA,\alpha)$ is an operator ideal, we let $\cA^f(X,Y)$ denote the
closure of $X^*\otimes Y$ under the norm $\alpha$.

In the formulas in this paper we shall, as is customary, interpret
$\pi_\infty$ as the operator norm and $i_\infty$ as the
$\gamma_\infty$-norm.

If $n\in\bN$ and $T\in B(\ell_2^n,X)$ then, following \cite{T}, we define
the $\ell$-norm of $T$ by
\[
\ell(T) = \left(\int_{\ell_2^n}\|Tx\|^2 d\gamma(x)\right)^{\frac12}
\]
where $\gamma$ is the canonical Gaussian probability measure on
$\ell_2^n$.

We let $(g_n)$ denote a sequence of independent standard Gaussian
variables on a fixed probability space $(\Omega,\cS,\sigma)$; it is
readily verified that if $T\in B(\ell_2^n,X)$ and $(\xi_j)$ denotes the
unit vector basis of $\ell_2$ then
\[
\ell(T) = \left(\int \|\sum^n_{j=1} g_j(t)T\xi_j\|^2
d\sigma(t)\right)^{\frac12}.
\]

A Banach space $X$ is said to have the Gordon-Lewis property
(abbreviated $\GL$) \cite{GL}, if every absolutely summing operator from
$X$ to an arbitrary Banach space $Y$ factors through $L_1$. It is
readily verified that $X$ has $\GL$ if and only if there is a constant
$K$ so that $\gamma_1(T)\le K \pi_1(T)$ for every Banach space $Y$ and
every $T\in X^*\otimes Y$. In that case $\GL(X)$ denotes the smallest
constant $K$ with this property.

We shall say that $X$ has $\GL_2$ if it has the above property with
$Y=\ell_2$ and we define the constant $\GL_2(X)$ correspondingly. An easy
trace duality argument yields that $\GL$ and $\GL_2$ are self dual
properties and that $\GL(X) = \GL(X^*)$, $\GL_2(X) = \GL_2(X^*)$ when
applicable. It is known that every Banach space with local unconditional
structure has $\GL$. For generalizations of $\GL$, see \cite{GJN}.

A Banach space $X$ is said to have the Gaussian Average property
(abbreviated $\GAP$) \cite{CN} if there is a constant $K$ so that
$\ell(T)\le K\pi_1(T^*)$ for every $T\in\ell_2^n\otimes X$ and every
$n\in\bN$. The smallest constant $K$ with this property is denoted
$\gap(X)$. 

A deep result of Pisier \cite{P2} states that a Banach space is
$K$-convex if and only if it is of type larger than 1. In this paper we
shall use this as the definition of $K$-convexity.

We shall also need some notation on operators with ranges in a Banach
lattice. Recall that if $Y$ is a Banach space and $L$ is a Banach
lattice then an operator $T\in B(Y,L)$ is called order bounded (see
e.g.\ \cite{S}, \cite{NJN1} and \cite{HNO}), if there exists a $z\in L$,
$z\ge 0$ so that
\begin{equation}
\label{eq0.1}
|Tx|\le \|x\|z\quad\mbox{for all $x\in Y$}
\end{equation}
and the order bounded norm $\|T\|_m$ is defined by
\begin{equation}
\label{eq0.2a}
\|T\|_m = \inf\{\|z\|\, \mid \, \mbox{$z$ can be used in
(\ref{eq0.1})}\}.
\end{equation}

It follows from \cite{K} and \cite{LT2} that if $T = \sum^n_{j=1}
y_j^*\otimes x_j\in Y^*\otimes L$ then
\[
\|T\|_m = \|\sup\{|\sum^m_{j=1} y^*_j(y)x_j| \, \mid \,\|y\|\le 1\}\|=\|
\, \|\sum^m_{j=1}x_jy^*_j\|_{Y^*}\|_L
\]
where the last equality is the definition of the 1-homogeneous
expression on the right.

We let $\cB(Y,L)$ denote the space of all order bounded operators from
$Y$ to $L$ equipped with the norm $\|\cdot\|_m$; it is readily seen to
be a Banach space and a left ideal.

If $X$ is a subspace of the Banach lattice $L$, then we let $Y\otimes_m
X$ denote the closure of $Y\otimes X$ in $\cB(Y^*,L)$ under the norm
$\|\cdot\|_m$. Note that $Y\otimes_m X$ depends on how $X$ is embedded
into $L$.

The next definition generalizes the concept of convexity and concavity
in Banach lattices.

\begin{definition}
\label{def0.1}
Let $X$ be a subspace of a Banach lattice $L$ and $1\le p< \infty$. $X$
is called $p$-convex in $L$ (respectively $p$-concave in $L$) if there
is a constant $K\ge 1$ so that for all finite sets
$\{x_1,x_2,\dots,x_n\}\subseteq X$ we have
\begin{equation}
\label{eq0.2}
\|(\sum^n_{j=1} |x_j|^p)^{\frac1p}\| \le K(\sum^n_{j=1}
\|x_j\|^p)^{\frac1p}
\end{equation}
(respectively
\begin{equation}
\label{eq0.3}
(\sum^n_{j=1} \|x_j\|^p)^{\frac1p} \le K\|(\sum^n_{j=1}
|x_j|^p)^{\frac1p}\|. )
\end{equation}
\end{definition}

The smallest constant $K$, which can be used in (\ref{eq0.2})
(respectively (\ref{eq0.3})) is denoted by $K^p(X,L)$ (respectively
$K_p(X,L)$). We put $K^p(L) = K^p(L,L)$ and $K_p(L) = K_p(L,L)$. Note
that the inequalities (\ref{eq0.2}) and (\ref{eq0.3}) depend on the
embedding of $X$ into $L$.

It follows from \cite{NJN1} that if $Y$ is a Banach space, $X$ is a
subspace of a Banach lattice $L$ and $T\in B(Y,X)$ with
$T^*\in\Pi_1^f(X^*,Y^*)$ then $T\in Y^*\otimes_m X$ with
$\|T\|_m\le\pi_1(T^*)$. The next theorem, which we shall use often in
the sequel generalizes this result (it also generalizes \cite[Theorem
1.3]{HNO} with an easier proof). Before we can state it we need a little
notation and a lemma.

Let $(\Delta,\cM,\mu)$ be a measure space, $X$ and $L$ as above and
$1\le p<\infty$. If $f\in L_p(\mu,X)$ is a simple function, say
$f=\sum^n_{j=1} 1_{A_j}x_j$, where $(x_j)^n_{j=1}\subseteq X$ and
$(A_j)^n_{j=1}$ is a set of mutually disjoint measurable sets then we
put
\[
\left(\int |f|^pd\mu\right)^{\frac1p} = \left(\sum^n_{j=1}\mu(A_j)
|x_j|^p\right)^{\frac1p}.
\]
The next lemma can be proved exactly as \cite[Proposition 1.2]{HNO}.

\begin{lemma}
\label{lemma0.2}
Let $1\le p<\infty$, $X$ a $p$-convex subspace of a Banach lattice $L$
and $f\in L_p(\mu,X)$. If $(s_n)\subseteq L_p(\mu,X)$ is a sequence of
simple functions with $\int \|f-s_n\|^p d\mu\to 0$ then $(\int |s_n|^p
d\mu)^{\frac1p}$ converges in $X$ to a limit, which only depends on $f$
and $p$. This limit is denoted by $(\int |f|^pd\mu)^{\frac1p}$ and
satisfies the inequalities:
\begin{eqnarray}
\label{eq0.4}
\|(\int |f|^p d\mu)^{\frac1p}\| &\le &
K^p(X,L)(\int\|f\|^pd\mu)^{\frac1p}\\
\label{eq0.5}
(\int|x^*(f)|^pd\mu)^{\frac1p} &\le & |x^*| ((\int|f|^pd\mu)^{\frac1p})
\quad\mbox{for all $x^*\in L^*$}.
\end{eqnarray}
\end{lemma}

We can now state

\begin{theorem}
\label{thm0.3}
Let $X$ be a $p$-convex subspace of a Banach lattice $L$, $1\le
p<\infty$, and let $Y$ be a Banach space. Then:
\begin{itemize}
\item[(i)]
\begin{equation}
\label{eq0.6}
\|T\|_m\le K^p(X,L)\pi_p(T^*)\quad\mbox{for all $T\in Y^*\otimes X$}.
\end{equation}
\item[(ii)] If $\, T\in B(Y,X)$ with $T^*\in\Pi_p(X^*,Y^*)$ then
$T\in\cB(Y,L^{**})$ with $\|T\|_m\le K^p(X,L)\pi_p(T^*)$.
\end{itemize}
\end{theorem}

\bproof
To prove (i) we let $T\in Y^*\otimes X$ and $\e>0$ be arbitrary. By
\cite[Lemma 1.8]{GJN} there is a finite dimensional subspace $F$ of $X$
containing $T(Y)$ so that if $T_F$ denotes $T$ considered as an operator
from $Y$ to $F$ then $\pi_p(T^*_F)\le \pi_p(T^*)+\e$.

By the Pietsch factorization theorem \cite{LT1} there exists a
probability measure $\mu$ on the unit ball $B_F$ of $F$ so that
\begin{equation}
\label{eq0.7}
\| T_F^*x^*\|\le\pi_p(T^*_F)(\int_{B_F}|x^*(x)|^p d\mu(x))^{\frac1p}.
\end{equation}
For every $y\in Y$ with $\|y\|\le 1$ and every $x^*\in L^*$, $x^*\ge 0$
we now get from (\ref{eq0.7}) and Lemma \ref{lemma0.2}:
\begin{equation}
\label{eq0.8}
|x^*(Ty)|\le \pi_p(T^*_F)(\int_{B_F}|x^*(x)|^p d\mu(x))^{\frac1p}\le
\pi_p(T^*_F)x^*((\int_{B_F}|x|^p d\mu(x))^{\frac1p}),
\end{equation}
which immediately gives
\begin{equation}
\label{eq0.9}
|Ty|\le\pi_p(T^*_F)(\int_{B_F}|x|^p d\mu(x))^{\frac1p}\quad\mbox{for all
$y\in B_Y$}.
\end{equation}
Hence
\begin{equation}
\label{eq0.10}
\|T\|_m\le\pi_p(T^*_F)\|(\int_{B_F}|x|^p d\mu(x))^{\frac1p}\|\le
K^p(X,L)(\pi_p(T^*))+\e)
\end{equation}
which gives (\ref{eq0.6}), since $\e$ was arbitrary.

(ii) can be proved in a similar manner. Noting that $X^{**}$ is
$p$-convex in $L^{**}$ with $K^p(X,L) = K^p(X^{**},L^{**})$ we get a
measure $\mu$ on $B_{X^{**}}$, so that
\[
|Ty|\le \pi_p(T^*)(\int_{B_{X^{**}}}|x^{**}|^p d\mu(x^{**}))^{\frac1p}
\]
where the right hand side represents an element in $L^{**}$.
\eproof

\section{Tensor Products and Embeddings of a Given Banach Space Into a
Banach Lattice}
\label{sec1}
\setcounter{equation}{0}

In this section we shall prove that every Banach space can be embedded
into a Banach lattice preserving a certain tensorial structure given a
priori. This result is based on an unpublished idea of Pisier \cite{P1},
but our construction is different
and is in nature similar to a result of Ruan \cite{Ru} on operator
spaces.

If $Y$ and $X$ are Banach spaces and $\alpha$ is a cross norm on
$Y\otimes X$ then we let $Y\otimes_\alpha X$ denote the completion of
$Y\otimes X$ under the norm $\alpha$.

If $E\subseteq Y$ and $F\subseteq X$ we can let $\alpha$ act on
$E\otimes F$ by considering it as a subspace of $Y\otimes X$ in the
canonical manner and define $E\otimes_\alpha F$ accordingly. Note
however that in general the outcome depends on how $E$, respectively
$F$, are embedded into $X$, respectively $Y$.

We make the following definition

\begin{definition}
\label{def1.1}
A crossnorm $\alpha$ on $Y\otimes X$ is called left tensorial if for all
$T\in B(Y)$ $T\otimes I_X\in B(Y\otimes X)$ with $\|T\otimes
I_X\|\le\|T\|$, where $I_X$ denotes the identity operator on $X$.
\end{definition}

\bremark
Note that the $m$-norm defined in section {\ref{sec0} is left tensorial.
\eremark

To obtain the main result of this section we shall be concerned with
left tensorial norms on $c_0\otimes X$ (or rather on
$\ell^n_\infty\otimes X$ for all $n\in\bN$). For technical reasons we
wish to have our left tensorial norms defined on $\ell_\infty\otimes X$
and hence need a few prerequisites. In passing we note that it is
fairly easy to see that if a norm $\alpha$ on $c_0\otimes X$ satisfies
the operator inequality in Definition \ref{def1.1} then it is a cross
norm up to a constant and hence left tensorial up to a constant.

In the rest of this section we let $(e_j)$ denote the unit vector basis
of $c_0$ with biorthogonal $(e_j^*)\subseteq\ell_1$; for all $n\in\bN$
$S^n_\infty$ denotes the unit sphere of $\ell^n_\infty$.

We need the following

\begin{proposition}
\label{prop1.2}
Let $X$ be a Banach space and $\alpha$ a left tensorial norm on
$c_0\otimes X$. If $u\in c_0\otimes_\alpha X$, then there is a unique
sequence $(x_n)\subseteq X$ so that
\begin{equation}
\label{eq1.1}
u = \sum^\infty_{n=1} e_n\otimes x_n
\end{equation}
and so that for all $n\in\bN$
\begin{equation}
\label{eq1.2}
\alpha\left(\sum^n_{j=1} e_j\otimes x_j\right)\le\alpha(u).
\end{equation}
\end{proposition}

\bproof
Let $P_n$ denote the natural projection of $c_0$ onto $[e_j\mid 1\le
j\le n]$ for every $n\in\bN$ and put $P_0=0$. By the left tensoriality
of $\alpha$ $P_n\otimes I_X$ is a bounded operator on $c_0\otimes X$ with
$\|P_n\otimes I_X\|=1$ so it admits an extension $Q_n\colon
c_0\otimes_\alpha X\to [e_j]^n_{j=1}\otimes X$ with $\|Q_n\|=1$ for all
$n\in\bN$. Note also that $\|Q_n-Q_{n-1}\|\le 1$. Let $u\in
c_0\otimes_\alpha X$ and put for every $n\in\bN$
\begin{equation}
\label{eq1.3}
e_n\otimes x_n=(Q_n-Q_{n-1})(u).
\end{equation}
For all $n\in\bN$ we have
\begin{equation}
\label{eq1.4}
Q_nu = \sum^n_{j=1}e_j\otimes x_j,
\end{equation}
from which (\ref{eq1.2}) follows.

Clearly $Q_nu\to u$ for all $u\in\Span\{e_j\}\otimes X$ and since
$\|Q_n\|=1$ for all $n\in\bN$ an easy density argument gives that
$Q_nu\to u$ for all $u\in c_0\otimes_\alpha X$ as well; this together
with (\ref{eq1.4}) gives (\ref{eq1.1}).
\eproof

>From this result we obtain:

\begin{proposition}
\label{prop1.3}
Let $X$ be a Banach space and $\alpha$ a left tensorial norm on
$c_0\otimes_\alpha X$. There exists a uniquely determined left tensorial
norm $\tilde{\alpha}$ on $\ell_\infty\otimes X$ so that
$\tilde{\alpha}|_{c_0\otimes X}=\alpha$. Here $c_0$ is considered as a
subspace of $\ell_\infty$ in the canonical manner.
\end{proposition}

\bproof
We consider the Banach space $(c_0\otimes_\alpha X)^{**}$ with its
canonical norm $\alpha^{**}$ and the idea is to identify $\ell_\infty
\otimes X$ with a canonical subspace of $(c_0\otimes_\alpha X)^{**}$ and
then put $\tilde{\alpha}$ equal to the restriction of $\alpha^{**}$ to
that subspace.

It is readily verified that $(c_0\otimes_\alpha X)^{*}$ can be
identified with the space $\ell_1\otimes_{\alpha^*} X^*$ consisting of
all sequences $(x_n^*)\subseteq X^*$ (written as $\sum_{n=1}^\infty
e_n^*\otimes x_n^*$) so that
\begin{equation}
\label{eq1.5}
\sum_{n=1}^\infty |x_n^*(x_n)|<\infty\quad\mbox{for all
$\sum^\infty_{n=1} e_n\otimes x_n\in c_0\otimes_\alpha X$}
\end{equation}
equipped with the norm
\begin{equation}
\label{eq1.6}
\alpha^*\left(\sum^\infty_{n=1} e_n^*\otimes x_n^*\right) = \sup\{
|\sum^\infty_{n=1} x_n^*(x_n)| \mid \alpha (\sum^\infty_{n=1} e_n\otimes
x_n)\le 1\}.
\end{equation}
Note that in particular we get for all $\sum^\infty_{n=1}
e_n^*\otimes x_n^*\in\ell_1\otimes_{\alpha^*}X^*$ and all $x\in X$:
\begin{equation}
\label{eq1.7}
\sum^\infty_{n=1} |x_n^*(x)|\le\alpha^*\left(\sum^\infty_{n=1}
e_n^*\otimes x_n^*\right) \|x\|
\end{equation}
and if $(\lambda_n)\in\ell_\infty$ with $|\lambda_n|\le 1$ for all
$n\in\bN$ then
$\sum^\infty_{n=1}e_n^*\otimes\lambda_nx^*_n\in\ell_1\otimes_{\alpha^*}X^*$
with
\begin{equation}
\label{eq1.8}
\alpha^*\left(\sum^\infty_{n=1} e_n^*\otimes\lambda_nx_n^*\right) =
\alpha^*\left(\sum^\infty_{n=1} e_n^*\otimes x_n^*\right).
\end{equation}
If $\sum^n_{j=1} h_j\otimes x_j\in\ell_\infty\otimes X$ and
$\sum^\infty_{i=1} e_i^*\otimes x_i^*\in\ell_1\otimes_{\alpha^*}X^*$
then by (\ref{eq1.7}) and (\ref{eq1.8})
\begin{equation}
\label{eq1.9}
\sum_{j=1}^n \sum_{i=1}^\infty |<h_j,e_i^*>| |x_i^*(x_j)|\le
\alpha^*\left(\sum_{i=1}^\infty e_i^*\otimes x_i^*\right)\sum_{j=1}^n \|x_j\|
\end{equation}
and hence we can let $\sum^n_{j=1} h_j\otimes x_j$ act as an element of
$(c_0\otimes_\alpha X)^{*}$ by the formula
\begin{equation}
\label{eq1.10}
\left<\sum_{j=1}^n h_j\otimes x_j,\sum_{i=1}^\infty e_i^*\otimes
x_i^*\right> = \sum_{j=1}^n\sum_{i=1}^\infty <h_j,e_i^*>
<x_i^*,x_j>\quad\mbox{for all $\sum_{i=1}^\infty e_i^*\otimes x_i^*
\in\ell_1\otimes_\alpha X$} 
\end{equation}
and we put $\tilde{\alpha}(\sum_{j=1}^n h_j\otimes x_j)$ equal to the
norm of that functional. Clearly $\tilde{\alpha}_{|c_0\otimes_\alpha
X}=\alpha$.

In order to prove that $\tilde{\alpha}$ is left tensorial we first note
that if $S\in B(c_0)$ then $(S\otimes I_X)^{**}=S^{**}\otimes I_X$ with
$\|S^{**}\otimes I_X\| = \|S\otimes I_X\| = \|S\|$. If $T\in B(\ell_\infty)$ is
arbitrary then since $\ell_\infty$ has the metric approximation
property, it follows from the local reflexitivity principle \cite{JRZ},
\cite{LR} that for every $\e>0$ there is a net $(S_t)$ of bounded
operators on $c_0$ with $\|S_t\|\le\|T\|+\e$ for all $t$ so that
\begin{equation}
\label{eq1.11}
\lim_t <S_t^{**}h,f> = <Th,f>\quad\mbox{for all $h\in\ell_\infty$ and
all $f\in\ell_1$}.
\end{equation}
For every $\sum_{i=1}^\infty e_i^* \otimes x_i^*
\in\ell_1 \otimes_{\alpha^*} X^*$ we now get:
\begin{eqnarray}
\label{eq1.12}
\lim_t \left\langle \sum_{j=1}^n S_t^{**} h_j \otimes x_j, \sum_{i=1}^\infty
e_i^* \otimes x_i^*\right\rangle &=& \sum_{j=1}^n \sum_{i=1}^{\infty} \lim_t <S_t^{**}
h_j,e_i^*> \nonumber \\
&=& \left\langle \sum_{j=1}^n (Th_j)\otimes x_j, \sum_{i=1}^\infty
e_i^* \otimes x_i^* \right\rangle
\end{eqnarray}
and hence $T\otimes I_X$ is bounded on
$\ell_\infty\otimes_{\tilde{\alpha}}X$ with
\begin{equation}
\label{eq1.13}
\|T\otimes I_X\|\le \|T\| +\e.
\end{equation}
It is clear that $\tilde{\alpha}$ is unique.
\eproof

The definition of left tensoriality immediately gives:

\begin{lemma}
\label{lemma1.4}
Let $X$ be a Banach space and $\alpha$ a left tensorial norm on
$\ell_\infty \otimes X$. If $F\subseteq \ell_\infty$ is a subspace and
$T\in B(F,\ell_\infty)$ then $T\otimes I_X\in B(F\otimes_\alpha
X,\ell_\infty \otimes_\alpha X)$ with $\|T\otimes I_X\| = \|T\|$. If $T$
is an isometry into then so is $T\otimes I_X$.
\end{lemma}

\bproof
Since $\ell_\infty$ has the extension property there is an extension
$\tilde{T}\colon \ell_\infty \to \ell_\infty$ of $T$ with
$\|\tilde{T}\|=\|T\|$ and hence the first part of the lemma follows from
the definition. If $T$ is an isometry into we can apply the same
procedure to $T^{-1}\colon T(F)\to F\subseteq \ell_\infty$ and the
result follows.
\eproof

We are now able to prove the finite dimensional version of our main
theorem:

\begin{theorem}
\label{thm1.5}
Let $E$ be an $n$-dimensional Banach space with a normalized basis
$(x_j)^n_{j=1}$ and biorthogonal system $(x_j^*)^n_{j=1}$ and let
$\alpha$ be a left tensorial norm on $\ell_\infty \otimes E$. There
exists a lattice norm $\|\cdot\|_\alpha$ on $C(S_\infty^n)$
\begin{equation}
\label{eq1.14}
ubc(x_j^*)^{-1} \|f\|_\infty \le \|f\|_\alpha \le
\alpha\left(\sum_{j=1}^n e_j\otimes x_j\right) \|f\|_\infty
\end{equation}
where $ubc(x_j^*)$ denotes the unconditional basis constant of
$(x^*_j)$, and an isometry $J\colon E\to (C(S_\infty^n),\|\cdot\|_\alpha)$ so that
for all $k\in\bN$ and all $y_1,y_2,\dots,y_k\in E$ we have
\begin{equation}
\label{eq1.15a}
\| \bigvee_{j=1}^k |Jy_j| \|_\alpha = \alpha\left(\sum_{j=1}^k
e_j\otimes y_j\right).
\end{equation}
\end{theorem}

\bproof
The construction of the norm $\|\cdot\|_\alpha$ is a kind of exercise
over the theme ``Krivine Calculus in Banach lattices'', \cite{K},
\cite{LT2}. We first note that if $u\in \ell_\infty \otimes E$, then
there are uniquely determined $h_1,h_2,\dots,h_n\in\ell_\infty$ so that
$u=\sum_{j=1}^n h_j\otimes x_j$. Further we let $\cS_n\subseteq
C(S_\infty^n)$ denote the set of all functions $p$ of the form
$p(t_1,t_2,\dots,t_n)=\|\sum_{j=1}^n t_jh_j\|_\infty$ for all
$(t_1,t_2,\dots,t_n)\in S^n_\infty$, where
$h_1,h_2,\dots,h_n\in\ell_\infty$. If $p\in\cS_n$ and
$h_1,h_2,\dots,h_n$ are as above then we shall say that
$(h_j)^n_{j=1}$represents $p$.

If $(f_j)^n_{j=1}\subseteq \ell_\infty$ also represents $p$, then
\begin{equation}
\label{eq1.15}
\| \sum^n_{j=1} t_jh_j\|_\infty = \|\sum_{j=1}^n
t_jf_j\|_\infty\quad\mbox{for all $(t_1,t_2,\dots,t_n)\in\bR^n$}
\end{equation}
and hence the operator $T\colon [h_j]\to[f_j]$ defined by
\[
T\left(\sum_{j=1}^n t_jh_j\right) = \sum^n_{j=1} t_jf_j\quad\mbox{for
all $t_1,t_2,\dots,t_n\in\bR$}
\]
is an isometry. From Lemma \ref{lemma1.4} it therefore follows that
\begin{equation}
\label{eq1.16}
\alpha\left(\sum_{j=1}^n h_j\otimes x_j\right) =
\alpha\left(\sum_{j=1}^n f_j\otimes x_j\right).
\end{equation}
Hence we can define $\|p\|_\alpha$ by
\begin{equation}
\label{eq1.17}
\|p\|_\alpha = \alpha\left(\sum_{j=1}^n h_j\otimes x_j\right).
\end{equation}
If $f\in C(S_\infty^n)$ then we define
\begin{equation}
\label{eq1.18}
\|f\|_\alpha = \inf\{\|p\|_\alpha \mid p\in\cS_n, |f|\le p\}.
\end{equation}
(Note that $|f(t_1,t_2,\dots,t_n)|\le \|f\|_\infty \|\sum_{j=1}^n
t_je_j\|_\infty$ for all $(t_1,t_2,\dots,t_n)\in S^n_\infty$).

>From Lemma \ref{lemma1.4} it follows in a similar manner as above that
if $p,q\in\cS_n$ with $p\le q$ then $\|p\|_\alpha\le\|q\|_\alpha$ and
therefore (\ref{eq1.18}) coincides with (\ref{eq1.17}) in case
$|f|\in\cS_n$. Thus $\|f\|_\alpha$ is well-defined for all $f\in
C(S^n_\infty)$.

We shall now show that $\|\cdot\|_\alpha$ is a norm on
$C(S^n_\infty)$. To this end let $p\in\cS_n$, $q\in\cS_n$ be represented
by $(h_j)^n_{j=1}\subseteq\ell_\infty$ and
$(f_j)^n_{j=1}\subseteq\ell_\infty$, respectively. We have to find a
representation of $p+q$. Put $F=([h_j]\oplus[f_j])_1$ and let $S$ be any
isometry of $F$ into $\ell_\infty$. For all
$(t_1,t_2,\dots,t_n)\in\cS_n$ we now get
\begin{eqnarray}
\label{eq1.19}
\|\sum^n_{j=1} t_j S(h_j,f_j)\|_\infty &=& \| S\left(\sum^n_{j=1}
t_jh_j,\sum^n_{j=1} t_jf_j\right)\|_\infty \nonumber \\
&=& \| \sum^n_{j=1} t_jh_j\|_\infty + \|\sum^n_{j=1}t_jf_j\|_\infty
\nonumber \\
&=& p(t_1,t_2,\dots,t_n)+q(t_1,t_2,\dots,t_n).
\end{eqnarray}
Hence using the definition of $\|\ \|_\alpha$ and Lemma \ref{lemma1.4}:
 \begin{eqnarray}
\label{eq1.20}
\|p+q\|_\alpha &=& \alpha\left(\sum^n_{j=1} S(h_j,f_j)\otimes
x_j\right)\nonumber \\
&= & \alpha\left(\sum^n_{j=1} S(h_j,0)\otimes x_j + \sum^n_{j=1}
S(0,f_j)\otimes x_j\right)\nonumber \\
&\le & \alpha\left(\sum^n_{j=1} S(h_j,0)\otimes
x_j\right)+\alpha\left(\sum^n_{j=1} S(0,f_j)\otimes x_j\right)\nonumber
\\
&=& \alpha\left(\sum^n_{j=1} h_j\otimes x_j\right) +
\alpha\left(\sum^n_{j=1} f_j\otimes x_j\right) = \|p\|_\alpha + \|q\|_\alpha.
\end{eqnarray}
Let now $f,g\in C(S^n_\alpha)$, let $\e>0$ be arbitrary and choose
$p,q\in\cS_n$ so that $|f|\le p$, $|g|\le q$, $\|p\|_\alpha\le
\|f\|_\alpha+\e$ and $\|q\|_\alpha\le\|g\|_\alpha+\e$. Since $|f+g|\le
p+q$ we obtain
\begin{equation}
\label{eq1.21}
\|f+g\|_\alpha\le \|p+q\|_\alpha\le \|p\|_\alpha+\|q\|_\alpha \le
\|f\|_\alpha + \|q\|_\alpha +2\e.
\end{equation}

Since $\e>0$ was arbitrary we have proved that $\|\cdot\|_\alpha$
satisfies the triangle inequality. It is clear that $\|af\|_\alpha = |a|
\|f\|_\alpha$ for all $f\in C(S^n_\infty)$ and all $a\in\bR$.

Let us now show the left inequality of (\ref{eq1.14}). Let $f\in
C(S^n_\infty)$ and $p\in\cS_n$ with $|f|\le p$. If
$(h_j)^n_{j=1}\subseteq\ell_\infty$ represents $p$ then we can define
$T\colon \ell_1\to E$ by $T = \sum^n_{j=1} h_j\otimes x_j$. For
arbitrary $(t_1,t_2,\dots,t_n)\in\cS^n_\infty$ we put $x^* =
\sum^n_{j=1} t_jx^*_j$ and get
\begin{eqnarray}
\label{eq1.22}
|f(t_1,t_2,\dots,t_n)| &\le & p(t_1,t_2,\dots,t_n) = \|\sum^n_{j=1}
x^*(x_j)h_j\|_\infty \nonumber \\
&=& \|T^*x^*\|_\infty \le \|T^*\| \| x^*\| \le \alpha(T)\|x^*\| =
\|p\|_\alpha \|\sum^n_{j=1} t_jx^*_j\|.
\end{eqnarray}
Taking first infimum over all $p\in\cS_n$ with $|f|\le p$ in
(\ref{eq1.22}) and thereafter supremum over all
$(t_1,t_2,\dots,t_n)\in\cS^n_\infty$ we obtain
\begin{equation}
\label{eq1.23}
\|f\|_\infty \le \|f\|_\alpha \sup\left\{\|\sum^n_{j=1} t_jx_j^*\|\mid
|t_j|=1\right\} = \|f\|_\alpha ubc(x_j^*).
\end{equation}
(\ref{eq1.23}) shows the left inequality of (\ref{eq1.14}) and hence
together with the above also gives that $\|\cdot\|_\alpha$ is a norm. It
follows immediately from the definition that if $f,g\in C(S^n_\infty)$
with $|f|\le |g|$ then $\|f\|_\alpha\le\|g\|_\alpha$ so that
$\|\cdot\|_\alpha$ is a lattice norm.

To prove the right inequality of (\ref{eq1.14}) we let again $f\in
C(S^n_\infty)$. Since for every $(t_1,t_2,\dots,t_n)\in S^n_\infty$ we
have
\begin{equation}
\label{eq1.24}
|f(t_1,t_2,\dots,t_n)|\le\|f\|_\infty \|\sum^n_{j=1} t_je_j\|_\infty
\end{equation}
we get by the definition of $\|f\|_\alpha$ that
\begin{equation}
\label{eq1.25}
\|f\|_\alpha \le \alpha\left(\sum^n_{j=1} e_j\otimes x_j\right) \|f\|_\infty
\end{equation}
which is the right inequality of (\ref{eq1.14}).

For every $1\le j\le n$ we let $\varphi_j\in C(S^n_\infty)$ be defined
by $\varphi_j(t_1,t_2,\dots,t_n)=t_j$ for all
$(t_1,t_2,\dots,t_n)\in\cS^n_\infty$ and define $J\colon E\to
C(S^n_\infty)$ by
\begin{equation}
\label{eq1.26}
J\left(\sum^n_{j=1} x_j^*(x)x_j\right)=\sum^n_{j=1}
x_j^*(x)\varphi_j\quad\mbox{for all $x\in E$}.
\end{equation}

We have to show that $J$ is an isometry and that (\ref{eq1.15}) holds.

To this end let $(h_j)^n_{j=1}\subseteq\ell_\infty$ and let $p\in\cS_n$
be represented by $(h_j)$. Since for every
$(t_1,t_2,\dots,t_n)\in\cS^n_\infty$ we have
\[
\sup \left\{ |\sum^n_{j=1} h^*(h_j)\varphi(t_1,t_2,\dots,t_n)| \mid
h^*\in\ell_1,\|h^*\|_1\le 1\right\} = \|\sum^n_{j=1} t_jh_j\|_\infty =
p(t_1,t_2,\dots,t_n)
\]
we get from the definition of $\|\cdot\|_\alpha$:
\begin{eqnarray}
\label{eq1.27}
\|\sum^n_{j=1} h_j\otimes J x_j\|_m &=& \|\sum^n_{j=1}
h_j\otimes\varphi_j\|_m\nonumber \\
&=& \|\sup_{\|h^*\|\le 1} |\sum^n_{j=1} h^*(h_j)\varphi_j| \|_\alpha
=\|p\|_\alpha = \alpha\left(\sum^n_{j=1} h_j\otimes x_j\right)
\end{eqnarray}
which is (\ref{eq1.15}) written in another form.

If $x\in E$ and we put $h_j=x_j^*(x)e_1$ for all $1\le j\le n$ in
(\ref{eq1.27}) we obtain that $J$ is an isometry.
\eproof

Before we can prove the main theorem of this section we need the
following proposition on the $m$-tensor product and ultraproducts of
Banach lattices.

\begin{proposition}
\label{prop1.6}
Let $(L_t)_{t\in I}$ denote a family of Banach lattices and let $L$
denote the Banach lattice obtained as the ultraproduct of $(L_t)$ along
an ultrafilter $\cU$.
For every $n\in\bN$ we have
\[
\ell^n_\infty\otimes_m L = \lim_\cU\ell^n_\infty\otimes_m L_t.
\]
\end{proposition}

\bproof
Let $Z=\{((x(t))\in\prod_{t\in I} L_t\mid \sup\{\|x(t)\| \mid t\in I\}
<\infty\}$ and let $\Phi\colon Z\to L$ denote the canonical quotient
map. Since by definition the ordering in $L$ is the one induced by
$\Phi$ it follows easily that if $n\in\bN$, $\{y_j\mid 1\le j\le
n\}\subseteq L$ and $(y_j(t))_{t\in I}\in Z$ with $\Phi((y_j(t))=y_j$
for all $1\le j\le n$ then $\Phi((\bigvee|y_j(t)|))=\bigvee^n_{j=1}|y_j|$, and
hence:
\begin{eqnarray}
\label{eq1.28}
\lim_\cU \|\sum^n_{j=1} e_j\otimes y_j(t)\|_m &=& \lim_\cU\|\bigvee^n_{j=1}
|y_j(t)| \|_{L_t}\nonumber \\
&=& \|\bigvee^n_{j=1} |y_j| \| = \|\sum^n_{j=1} e_j\otimes y_j\|_m.
\end{eqnarray}
\eproof

We can now easily prove

\begin{theorem}
\label{thm1.7}
Let $X$ be a Banach space and $\alpha$ a left tensorial norm on
$\ell_\infty \otimes X$. There exist a Banach lattice $L$ and an
isometry $J$ of $X$ into $L$ so that for all $k\in\bN$ and all
$y_1,y_2,\dots,y_k\in X$ we have
\begin{equation}
\label{eq1.29}
\|\bigvee^k_{j=1} |J y_k| \| = \alpha\left(\sum^k_{j=1} e_j\otimes y_j\right).
\end{equation}
\end{theorem}

\bproof
For every finite dimensional subspace $E\subseteq X$ we consider
$\ell_\infty\otimes E$ as a subspace of $\ell_\infty \otimes X$ in the
natural way and equip it with the norm $\alpha$ restricted to
$\ell_\infty \otimes E$, i.e.\ we put $\ell_\infty \otimes_\alpha E =
(\ell_\infty \otimes E,\alpha)$ (this is a slight misuse of notation
which can cause problems for concrete $\alpha$'s but we shall only use
it in this proof). $\alpha$ is clearly left tensorial on $\ell_\infty
\otimes E$.

Put
\begin{equation}
\label{eq1.30}
\cF = \{E\subseteq X\mid E \ \mbox{finite dimensional$\}$}.
\end{equation}
In every $E\in\cF$ we choose a normalized basis and let $L_E$ be the
Banach lattice constructed in Theorem \ref{thm1.5} relative to the
chosen basis and our choice of $\ell_\infty \otimes_\alpha E$ and let
$J_E\colon E\to L_E$ be the isometry constructed there.

We define $L$ to be the ultraproduct of $\{L_E\mid E\in \cF\}$ along a
free ultrafilter $\cU$ of $\cF$.

Let $Q$ be the canonical quotient map of $(\Pi L_E)_\infty$ onto $L$ and
let for every $x\in E$ $\tilde{x}\in(\Pi L_E)_\infty$ be defined by
\begin{equation}
\label{eq1.31}
\tilde{x}(E) = \left\{\begin{array}{cc} J_Ex & \mbox{if $x\in E$}\\0 &
\mbox{else} \end{array}\right. \quad\mbox{for every $E\in\cF$}.
\end{equation}
If we put $Jx = Q\tilde{x}$ for all $x\in X$, then $J$ is readily seen
to be a linear map from $X$ to $L$ and if $x\in E$ then it follows from
(\ref{eq1.31}) and the definition of the norm in $L$ that
\begin{equation}
\label{eq1.32a}
\|Jx\| = \lim_\cU \|\tilde{x}(E)\| = \|x\|.
\end{equation}

Hence $J$ is an isometry of $X$ into $L$.

If finally $n\in\bN$ and $\{x_j\mid 1\le j\le n\}\subseteq X$ then it
follows from Theorem \ref{thm1.5} and Proposition \ref{prop1.6} that
\begin{equation}
\label{eq1.32}
\|\bigvee^n_{j=1} |J x_j| \| = \lim_\cU \|\bigvee^n_{j=1}
\tilde{x}_j(E)\|_{L_E} = \alpha\left(\sum^n_{j=1} e_j\otimes x_j\right).
\end{equation}
(Note that our special choice of $\ell_\infty \otimes_\alpha E$ is
important here.)
\eproof

We end this section with a few corollaries:

\begin{corollary}
\label{cor1.8}
If $X$ and $L$ are as in Theorem \ref{thm1.7} then for all $n\in\bN$,
all $\{x_j\mid 1\le j\le n\}\subseteq X$ and all $\{h_j \mid 1\le j\le
n\}\subseteq\ell_\infty$ we have
\begin{equation}
\label{eq1.33}
\|\sum^n_{j=1} h_j\otimes Jx_j\|_m = \alpha\left(\sum^n_{j=1} h_j\otimes
x_j\right). 
\end{equation}
\end{corollary}

\bproof
Let $n\in\bN$, $\{h_j\mid 1\le j\le n\}\subseteq\ell_\infty$, $\{x_j\mid
1\le j\le n\}\subseteq X$ and $\e>0$ be given. Using the local
properties of $\ell_\infty$ we can find an $m\in\bN$ and an isomorphism
$T$ of $[h_j]$ into $\ell_\infty^m$ so that $\|T\|=1$ and $\|T^{-1}\|\le
(1+\e)$.

By Lemma \ref{lemma1.4} we get:
\begin{equation}
\label{eq1.34}
(1+\e)^{-1} \alpha\left(\sum^n_{j=1} h_j\otimes x_j\right) \le
\alpha\left(\sum^n_{j=1} T h_j\otimes x_j\right)\le
\alpha\left(\sum^n_{j=1} h_j\otimes x_j\right)
\end{equation}

\begin{equation}
\label{eq1.34a}
(1+\e)^{-1} \left\|\sum^n_{j=1} h_j\otimes J
x_j\right\|_m\le\left\|\sum^n_{j=1} T h_j\otimes J
x_j\right\|_m\le\left\|\sum^n_{j=1} h_j\otimes J x_j\right\|_m.
\end{equation}

Since by Theorem \ref{thm1.7} $\|\sum^n_{j=1} T h_j\otimes J
x_j\|_m=\alpha(\sum^n_{j=1} T h_j\otimes x_j)$ we get (\ref{eq1.33}) by
letting $\e$ tend to 0.
\eproof

The next corollary follows from trace class duality.

\begin{corollary}
\label{cor1.9}
Let $X$ and $L$ be as in Theorem \ref{thm1.7} and let $\alpha^*$ be the dual
norm to $\alpha$. Every operator $T\colon J(X)\to\ell_1$ with
$TJ\in(\ell_\infty\otimes_\alpha X)^*$ extends to an operator
$\widetilde{T}\colon L\to\ell_1$ with
$\widetilde{T}~^*\in\cB(\ell_\infty,L^*)$ so that
\begin{equation}
\label{eq1.38}
\|T^*\|_m = \alpha^*(TJ).
\end{equation}
\end{corollary}

\bproof
>From \cite{S} it follows that
$(\ell_\infty\otimes_mL)^*=\cB(\ell_\infty,L^*)$ and since
$\ell_\infty\otimes_\alpha X$ is canonically isometric to
$\ell_\infty\otimes_m J(X)$ we get that the restriction map from $L^*$
onto $J(X)^*$ induces a quotient map of $\cB(\ell_\infty,L^*)$ onto
$(\ell_\infty\otimes_\alpha X)^*$. Hence (\ref{eq1.38}) follows.
\eproof

\bremark
Using the local properties of $L_1$-spaces it is readily verified that
Corollary \ref{cor1.9} still holds if $\ell_1$ is substituted by
$L_1(\mu)$, where $\mu$ is an arbitrary measure.
\eremark

\section{Some Applications}
\label{sec2}
\setcounter{equation}{0}

In this section we shall give some applications of Theorem \ref{thm1.7}
and its corollaries. We start with

\begin{theorem}
\label{thm2.1}
Let $X$ be a Banach space. Then there exists a Banach space $L$ so that
$X$ embeds isometrically into $L$ (we write $X\subseteq L$) and so that
for all Banach spaces $F$ we have
\begin{equation}
\label{eq2.1}
T\in F \otimes_m X \Longleftrightarrow T^*\in\Pi_1^f(X^*,F).
\end{equation}
\end{theorem}

\bproof
Put $\ell_\infty \otimes_\alpha X = \ell_\infty \otimes_\pi X$ and let
$L$ be the Banach lattice constructed in Theorem
\ref{thm1.7} so that $X\subseteq L$ and $\ell_\infty \otimes_m X = \ell_\infty \otimes_\pi
X$. To prove that $L$ has the desired property it is enough to prove
(\ref{eq2.1}) when $F\subseteq\ell_\infty$. We clearly have
\[
T\in F\otimes_m X\Longleftrightarrow (T^*\colon X^*\to \ell_\infty)\in
N_1(X^*,\ell_\infty)\Longleftrightarrow T^*\in\Pi_1^f(X^*,F)
\]
which is (\ref{eq2.1}).
\eproof

In analogy with Corollary \ref{cor1.9} we get the following corollary,
using trace class duality arguments.

\begin{corollary}
\label{cor2.2}
Let $X$ be a Banach space and let $L$ be the Banach lattice constructed
in Theorem \ref{thm2.1}. If $G$ is another Banach space then every
operator $T\in\Gamma_1(X,G)$ admits an extension
$\tilde{T}\in\Gamma_1(L,G)$ with
$\gamma_1(T)=\gamma_1(\tilde{T})$. Furthermore,
\begin{equation}
\label{eq2.2}
\gamma_1(\tilde{T}) = \inf\|A\| \|B\| \mid \exists \ \mbox{a measure}\ \nu,
A\colon L\to L_1(\nu), A\ge 0, B\colon L_1(\nu)\to G, T=BA\}.
\end{equation}
\end{corollary}

\bproof
Let $T\in\Gamma_1(X,G)$, let $\e>0$ be arbitrary and choose a measure
$\mu$ and operators $S\colon X\to L_1(\mu)$, $U\colon L_1(\mu)\to G$ so
that $T=US$, $\|U\|\le 1$ and $\|S\|\le \gamma_1(T)+\e$. By Corollary
\ref{cor1.9} and its remark $S$ admits an extension $\tilde{S}\colon
L\to L_1(\mu)$ so that $\tilde{S}^*\in\cB(L_\infty(\mu),L^*)$ and
$\|S\|=\|\tilde{S}\|_m$. Since $X^*$ is order complete and $L_1(\mu)$ is
complemented in $L_1(\mu)^{**}$ it follows e.g. from \cite{NJN1} that there exists a measure $\nu$, a positive operator $A\colon
L\to L_1(\nu)$ and an operator $V\colon L_1(\nu)\to L_1(\mu)$ so that
$\tilde{S}=VA$ and
\begin{equation}
\label{eq2.3}
\|A\| \|V\| \le \|\tilde{S}^*\|_m +\e\le \|S\|+\e.
\end{equation}
The operator $\tilde{T}=U\tilde{S}$ clearly extends $T$ and belongs to
$\Gamma_1(L,G)$. Furthermore, $\tilde{T}=UVA$ and hence
\begin{eqnarray}
\label{eq2.4}
\gamma_1(T) &\le & \gamma_1(\tilde{T})\le \|A\| \|UV\| \le \|A\| \|U\|
\|V\| \nonumber \\
& \le  & \|S\| \|U\| + \e\|U\|\le \gamma_1(T)+2\e
\end{eqnarray}
so that $\gamma_1(T)=\gamma_1(\tilde{T})$ and (\ref{eq2.2}) holds.
\eproof

The next theorem characterizes Banach spaces with $\GL_2$ in terms of
embeddings into Banach lattices. It generalizes \cite[Corollary
2.3]{CN}.

\begin{theorem}
\label{thm2.3}
Let $X$ be a Banach space. The following two statements are equivalent
\begin{itemize}
\item[(i)] $X$ has $\GL_2$.
\item[(ii)] There exists a Banach lattice $L\supseteq X$ so that every
$T\in\Pi_1(X,\ell_2)$ admits an extension $\tilde{T}\in\Pi_1(L,\ell_2)$.
\end{itemize}
\end{theorem}

\bproof
Since every Banach lattice has $\GL_2$ by \cite{GL} (ii) trivially 
implies (i). Next, assume that $X$ has $\GL_2$ and let $L$ be the Banach
lattice constructed in Theorem \ref{thm2.1}. If $T\in\Pi_1(X,\ell_2)$
then also $T\in\Gamma_1(X,\ell_2)$ with $\gamma_1(T)\le\GL_2(X)\pi_1(T)$
and hence by Corollary \ref{cor2.2} $T$ admits an extension
$\tilde{T}\in\Gamma_1(L,\ell_2)$ with
$\gamma_1(\tilde{T})=\gamma_1(T)$. However, by Grothendieck's theorem
\cite{LT1}, $\tilde{T}$ is also 1-summing with
\begin{equation}
\label{eq2.5}
\pi_1(\tilde{T}) \le K_G \gamma_1(\tilde{T})\le K_G \GL_2(X)\pi_1(T)
\end{equation}
where $K_G$ is Grothendieck's constant.
\eproof

\bremark
>From Theorem \ref{thm2.1} and Corollary \ref{cor2.2} it follows that if
$X$ has $\GL_2$ and $T\in\Pi_1(X,\ell_2)$ then $\gamma_1(T)$ can be
computed by looking on factorizations $T=BA$ where $A$ is the
restriction to $X$ of a positive operator from a suitable Banach lattice
$L\supseteq X$ to an $L_1$-space.

Let us note in passing that if $X$ is contained in a Banach lattice
$L$ and (ii) of Theorem \ref{thm2.3} holds then
\[
\ell_2\otimes_m X = \{ T\colon\ell_2\to X\mid
T^*\in\Pi^f_1(X^*,\ell_2)\}.
\]
Indeed, it easily follows that there is a constant $K\ge 1$ so that
every $S\in\Pi_1(X,\ell_2)$ admits an extension
$\tilde{S}\in\Pi_1(L,\ell_2)$ with $\pi_1(\tilde{S})\le K
\pi_1(S)$. Hence if $S\in\Gamma_1(X,\ell_2)$, $S$ admits an extension
$\tilde{S}\in\Gamma_1(L,\ell_2)$ with
$\gamma_1(\tilde{S})\le\pi_1(\tilde{S})\le K\pi_1(S)\le
KK_G\gamma_1(S)$. An easy trace duality argument now gives that if $Q$
denotes the natural quotient map of $L^*$ onto $X^*$, then for every
$T\colon\ell_2\to X$ with $T^*\in\Pi_1(X^*,\ell_2)$ we have
$\pi_1(T^*)\le KK_G\pi_1(T^*Q)$.

If now $T\in\ell_2\otimes_m X$ then by \cite{HNO} and the above
$\pi_1(T^*)\le KK_G\pi_1(T^*Q)\le KK_G^2\|T\|_m$, so that
$T^*\in\Pi^f_1(X^*,\ell_2)$. The other direction follows from Section
\ref{sec0}, since if $T^*\in\Pi^f_1(X^*,\ell_2)$, $T\in\ell_2\otimes_mX$
with $\|T\|_m\le \pi_1(T^*Q)\le \pi_1(T^*)$.
\eremark

We can now characterize Banach spaces with $\GAP$.

\begin{theorem}
\label{thm2.4}
Let $X$ be a Banach space. The following statements are equivalent:
\begin{itemize}
\item[(i)] $X$ has $\GAP$.
\item[(ii)] There exist a Banach lattice $L\supseteq X$ and a constant
$K\ge 1$ so that for all $x_1,x_2,\dots,x_n\in X$ we have
\begin{equation}
\label{eq2.6}
\|\left(\sum^n_{j=1} |x_j|^2\right)^{\frac12}\|\le
\left(\int\|\sum^n_{j=1} g_i(t)x_i\|^2 d\sigma(t)\right)^{\frac12} \le
K\|\left(\sum^n_{j=1} |x_i|^2\right)^{\frac12}\|.
\end{equation}
\end{itemize}
\end{theorem}

\bproof
(i)  $\Rightarrow$ (ii). Let $L$ be the Banach lattice constructed in
Theorem \ref{thm2.1}, and let $(f_j)$ denote the unit vector basis of
$\ell_2$. If $\{x_1,x_2,\dots,x_n\}\subseteq X$ and $T=\sum^n_{j=1}
f_j\otimes x_j$ then it follows from Theorem \ref{thm2.1} and the GAP of
$X$ that
\begin{eqnarray}
\label{eq2.7}
\left(\int\|\sum^n_{j=1} g_j(t)x_j\|^2 d\sigma(t)\right)^{\frac12} &=&
\ell(T) \le \gap(X)\pi_1(T^*)\nonumber \\
&\le & \gap(X)\|T\|_m = \gap(X)\|\left(\sum^n_{j=1} |x_j|^2\right)^{\frac12}\|.
\end{eqnarray}
The left inequality of (\ref{eq2.6}) always holds in a Banach lattice
\cite{LT2}. This shows (i) $\Rightarrow$ (ii).

(ii) $\Rightarrow$ (i). Assume that $X\subseteq L$ and that
(\ref{eq2.6}) holds. If $T\in\ell_2\otimes X$ then
\begin{equation}
\label{eq2.8}
\ell(T) \le K\|T\|_m \le K\pi_1(T^*)
\end{equation}
so that $X$ has $\GAP$.
\eproof

\bremark
It follows from \cite{CN} that a Banach space with $\GAP$ is of finite
cotype and one could hope that the Banach lattice in Theorem
\ref{thm2.4} could be chosen to be of finite cotype. However this is not
the case. Indeed, \cite[Example 1.16]{CN} shows that the Schatten class
$c_p$ for $2<p<\infty$ has GAP but not $(S)$ and therefore a Banach
lattice $L\supseteq c_p$ with the properties of Theorem \ref{thm2.4}
cannot be of finite cotype, since every subspace of a Banach lattice of
finite cotype has $(S)$.

We also note that if a Banach space $X$ is contained in a Banach lattice
$L$, so that (\ref{eq2.6}) holds then it follows from \cite[Proposition
0.3]{CN} that there is a constant $K_1\ge 1$ so that for all
$T\in\ell_2\otimes X$ we have $K_1^{-1} \pi_1(T^*)\le \ell(T)\le
\|T\|_m\le K\pi_1(T^*)$. Furthermore, an easy trace duality argument,
similar to the one in Corollary \ref{cor2.2}, applied to these
inequalities yields that every operator $T\in\Gamma_1(X,\ell_2)$ admits
an extension $\tilde{T}\in\Gamma_1(L,\ell_2)$.  
\eremark

Combining Theorem \ref{thm2.4} with \cite[Theorem 1.9]{CN} we obtain

\begin{corollary}
\label{cor2.5}
Let $X$ be a Banach space. The following statements are equivalent
\begin{itemize}
\item[(i)] $X$ is $K$-convex and has $\GL_2$.
\item[(ii)] There exist Banach lattices $L\supseteq X$, $Y\supseteq X^*$
and a constant $K\ge 1$ so that the inequality (\ref{eq2.6}) holds for
finite sets of vectors in $X$, respectively in $X^*$.
\end{itemize}
\end{corollary}

It follows from the remark above that if (ii) of Corollary \ref{cor2.5}
holds then every operator $T\in\Pi_1(X,\ell_2)$ admits an extension
$\tilde{T}\in\Pi_1(L,\ell_2)$. Note also that since $X$ is $K$-convex in
that case (\ref{eq2.6}) shows that $\ell_2\otimes_m X^*$ is canonically
isomorphic to $(\ell_2\otimes_m X)^*$.

We now introduce a property of Banach spaces which is more general than
property $(S)$ defined in \cite{CN}.

\begin{definition}
\label{def2.6}
Let $1\le p\le q<\infty$. A Banach space $X$ is said to have $S\GL(p,q)$
if there is a constant $K$ so that if $Y$ is an arbitrary Banach space,
then
\begin{equation}
\label{eq2.9}
\pi_q(T)\le K \pi_p(T^*)\quad\mbox{for all $T\in Y^*\otimes X$}.
\end{equation}
\end{definition}

If we put $Y=\ell_2$ and $p=1$ in this definition we get the property
$(S_q)$ of \cite{CN}. It is immediate that subspaces of Banach spaces
with property $\GL(p,q)$ from \cite{GJN} have $S\GL(p,q)$. In particular it follows from
\cite[Theorem 1.3]{GJN} (see also \cite{HNO} and \cite{J}) that every
subspace of a $p$-convex and $q$-concave Banach lattice has
$S\GL(p,q)$. It is actually also a consequence of our next result.

We now wish to characterize property $S\GL(p,q)$ in terms of embeddings
into Banach lattices. The result states:

\begin{theorem}
\label{thm2.7}
Let $X$ be a Banach space and $1\le p\le q<\infty$. The following
statements are equivalent:
\begin{itemize}
\item[(i)] $X$ has $S\GL(p,q)$.
\item[(ii)] $X$ satisfies (\ref{eq2.9}) with $Y=\ell_{q'}$.
\item[(iii)] There exists a Banach lattice $L$ with $X\subseteq L$ so
that $X$ is $p$-convex and $q$-concave in $L$.
\end{itemize}
\end{theorem}

\bproof
(i) $\Rightarrow$ (ii) is obvious so assume that (ii) holds and let $L$
be the Banach lattice constructed in Theorem \ref{thm1.7} with
\begin{equation}
\label{eq2.10}
T\in \ell_\infty \otimes_\alpha X\Longleftrightarrow
T^*\in\Pi^f_p(X^*,\ell_\infty)
\end{equation}
as the defining tensor norm.

Further, let $(f_j)$ denote the unit vector basis of $\ell_q$ with
biorthogonal system $(f^*_j)\subseteq\ell_{q'}$. If $n\in\bN$ and
$\{x_1,x_2,\dots,x_n)\subseteq X$ then with $T=\sum^n_{j=1} f_j\otimes
x_j$ we obtain:
\begin{eqnarray}
\label{eq2.11}
\left(\sum^n_{j=1} \|x_j\|^q\right)^{\frac1q} &=& \left(\sum^n_{j=1}
\|Tf^*_j\|^q\right)^{\frac1q}\nonumber \\
&\le& \pi_q(T)\le K\pi_p(T^*)=K\|T\|_m=K\|\left(\sum^n_{j=1}
|x_j|^q\right)^{\frac1q}\|,
\end{eqnarray}
which shows that $X$ is $q$-concave in $L$.

If $(u_i)$ denotes the unit vector basis in $\ell_p$ we get with 
$U=\sum u_i\otimes x_i$
\begin{equation}
\label{eq2.12}
\|\left(\sum^n_{j=1} |x_j|^p\right)^{\frac1p}\| = \|U\|_m =
\pi_p(U^*)\le\left(\sum^n_{j=1} \|x_j\|^p\right)^{\frac1p},
\end{equation}
which shows that $X$ is $p$-convex in $L$.

Assume next that (iii) holds, put $K=K_q(X,L)$ and let $Y$ be an
arbitrary Banach space, $T\in Y^*\otimes X$ and $\e>0$.

Since $L^{**}$ is order complete there is a compact Hausdorff space
$\Delta$ and operators \linebreak $A\in B(Y,C(\Delta))$, $B\in
B(C(\Delta),L^{**})$, $B\ge 0$ so that $\|A\| \|B\|\le \|T\|_m+\e$ and
$T=BA$.

If $y_1,y_2,\dots,y_n\in Y$ are arbitrarily chosen, then since $B\ge 0$:
\begin{eqnarray}
\label{eq2.13}
\left(\sum^n_{j=1} \|Ty_j\|^q\right)^{\frac1q} &\le &
K\|\left(\sum^n_{j=1} |BAy_j|^q\right)^{\frac1q}\| \nonumber \\
&\le & K\|B\| \|\left(\sum^n_{j=1} |A(y_j)|^q\right)^{\frac1q}\| = K\|B\|
\sup_{t\in \Delta} \left(\sum^n_{j=1} |(Ay_j)(t)|^q\right)^{\frac1q}\nonumber
\\
&=& K\|B\| \sup\left\{\left(\sum^n_{j=1}
|\mu(Ay_j)|^q\right)^{\frac1q}\mid \mu\in C(\Delta)^*,\|\mu\|\le
1\right\}\nonumber \\
&\le & K\|A\| \|B\| \sup\left\{\left(\sum^n_{j=1}
|y^*(y_j)|^q\right)^{\frac1q}\mid y^*\in Y^*,\|y^*\|\le 1\right\}.
\end{eqnarray}
This shows that $\pi_q(T)\le K\|A\| \|B\| \le K(\|T\|_m+\e)$ and hence
since $\e>0$ was arbitrary
\begin{equation}
\label{eq2.14}
\pi_q(T)\le K\|T\|_m.
\end{equation}
Furthermore, from Theorem \ref{thm0.3} it follows that $\|T\|_m\le
K^p(X,L)\pi_p(T^*)$; thus concluding that $X$ has $S\GL(p,q)$.
\eproof

We end this section by giving a characterization of the $\GL$-property
in terms of Banach lattices, but it is less intuitive than the results
above. Before we can formulate it we need a little notation.

If $X$ is a subspace of a Banach lattice $L$ and $E$ is a finite
dimensional Banach space we denote the norm in $(E\otimes_m X)^*$ by
$\|\cdot\|^*_m$. If $S\in E^*\otimes X^*$ and $T\in E\otimes_m X$ then
\begin{equation}
\label{eq2.15}
|\Tr(T^*S)|\le \nu_1(T^*S)\le \nu_1(S)\|T\|_m
\end{equation}
and therefore $S$ acts as a bounded linear functional on $E\otimes_m X$
by the formula
\begin{equation}
\label{eq2.16}
\langle S,T\rangle = \Tr(T^*S)\quad\mbox{for all $T\in E\otimes_mX$}
\end{equation}
and
\begin{equation}
\label{eq2.17}
\|S\|^*_m = \sup\{|\Tr(T^*S)| \mid T\in E\otimes_m X, \|T\|_m\le 1\}.
\end{equation}
We are now able to prove the following:

\begin{theorem}
\label{thm2.8}
Let $X$ be a Banach space. The following statements are equivalent:
\begin{itemize}
\item[(i)] $X$ has $\GL$.
\item[(ii)] There exist Banach lattices $L$ and $M$ so that $X\subseteq
L$, $X^*\subseteq M$ and a constant $K\ge 1$ so that for all finite
dimensional Banach spaces $E$ we have
\begin{equation}
\label{eq2.18}
\| S \|^*_m\le K\|S\|_m\quad\mbox{for all $S\in E^*\otimes X^*$}.
\end{equation}
\end{itemize}
\end{theorem}

\bproof
Assume first that $X$ has $\GL$ with $\GL$-constant $K$, and let $L$ and
$M$ be the Banach lattices constructed in Theorem \ref{thm2.1} with
$X\subseteq L$ and $X^*\subseteq M$. If $E$ is an arbitrary finite
dimensional Banach space and $T\in E\otimes_m X$ then by the definition
of $L$ and the $\GL$-property we have
\begin{equation}
\label{eq2.19}
\gamma_\infty(T) = \gamma_1(T^*)\le K\pi_1(T^*)=K\|T\|_m.
\end{equation}
Hence by duality and the definition of $M$ we get for every $S\in
E^*\otimes X^*$:
\begin{equation}
\label{eq2.20}
\|S\|^*_m\le K\pi_1(S^*)=K\|S\|_m,
\end{equation}
which shows that the inequality in (ii) holds.

Assume next that (ii) holds and let $E$ be an arbitrary finite
dimensional Banach space. Since for every $T\in E^*\otimes X$ we have
$\|T\|_m\le \pi_1(T^*)$ it follows by duality that for all $S\in
E\otimes X^*$ we have $\gamma_\infty(S)\le\|S\|^*_m$. Hence by
(\ref{eq2.18}) we get for every $T\in X^*\otimes E$
\begin{equation}
\label{eq2.21}
\gamma_1(T) = \gamma_\infty(T^*)\le \|T^*\|_m^* \le K\|T^*\|_m\le K\pi_1(T)
\end{equation}
which shows that $X$ has $\GL$ with constant less than or equal to $K$.
\eproof

\bremark
By putting $E=\ell_2$ in Theorem \ref{thm2.8} it is readily verified
that a similar result holds for $\GL_2$. In (ii) we actually get that
the norms $\|\cdot\|_m^*$ and $\|\cdot\|_m$ become equivalent. The $\GL_2$-version of
Theorem \ref{thm2.8} can also easily be derived from Theorem \ref{thm2.3}.
\eremark

Recall that an operator $T$ from a Banach lattice $L$ to a Banach space
$Y$ is called cone-summing if it maps unconditional convergent series 
of positive vectors to absolutely convergent ones. It follows from \cite{S}
that $T$ is cone-summing if and only if $T^*$ is order bounded. Using the
same argument as in the proof of Theorem \ref{thm2.3} we easily obtain:

\begin{theorem}
\label{thm2.9}
If $X$ is a Banach space, then the following two statements are
equivalent:
\begin{itemize}
\item[(i)] $X$ has $\GL$.
\item[(ii)] There exists a Banach lattice $L\supseteq X$ so that every 
absolutely summing operator $T$ from $X$ to an arbitrary Banach space
$Y$ extends to a cone-summing operator $\tilde{T}$ from $L$ to $Y$.
\end{itemize}
\end{theorem}

\section{$\GL$-subspaces of Banach lattices of finite cotype}
\label{sec3}

\setcounter{equation}{0}

As noted in the previous section it is not always possible to embed a
Banach space of finite cotype into a Banach lattice of finite cotype
(equivalently of finite concavity), not even if it has GAP. However,
combining the results of this paper with those of \cite{CN} we believe
that the following two conjectures have positive answers.

\begin{conjecture}
\label{conj3.1}
If a Banach space $X$ has $S\GL(1,q)$ for some $q$, $1\le q<\infty$,
then $X$ embeds into a $q$-concave Banach lattice $L$.
\end{conjecture}

\begin{conjecture}
\label{conj3.2}
A $K$-convex Banach space $X$ has $\GL_2$ if and only if both $X$ and
$X^*$ embed into Banach lattices of finite cotype.
\end{conjecture}

In this section we shall investigate when a subspace of a Banach lattice
has $\GL$ or $\GL_2$. For convenience we shall say that a subspace $X$
of a Banach lattice $L$ is optimally embedded into $L$ if
$\ell_2\otimes_m X = \Gamma_\infty^f(\ell_2,X)$.

It follows from Theorem \ref{thm2.3} and the remark just after that a
Banach space $X$ has $\GL_2$ if and only if there exists a Banach
lattice $L$ so that $X$ can be optimally embedded into $L$. On the other
hand it follows from \cite[Corollary 2.3]{CN} that if $X$ is
a subspace of a Banach lattice of finite cotype then $X$ has $\GL_2$ if
and only if it is optimally embedded into $L$; in that case any other
embedding of $X$ into $L$ is also optimal.

>From the results in Section \ref{sec2} we can conclude

\begin{theorem}
\label{thm3.2}
Let $X$ be a subspace of a Banach lattice $L$ of finite concavity. The
following statements are equivalent:
\begin{itemize}
\item[(i)] $X$ is $K$-convex and has $\GL_2$.
\item[(ii)] There exist a Banach lattice $M$ with $X^*\subseteq M$ and a
constant $K\ge 1$ so that
\begin{eqnarray}
\label{eq3.1}
K^{-1}\|\left(\sum^n_{j=1} |x_j^*|^2\right)^{\frac12}\| &\le &
\left(\int\|\sum^n_{j=1} g_j(t)x^*_j\|^2 d\sigma(t)\right)^{\frac12}
\nonumber \\
&\le & K\|\left(\sum^n_{j=1} |x^*_j|^2\right)^{\frac12}\|
\end{eqnarray}
for all finite sets $\{x^*_j\mid 1\le j\le n\}\subseteq X^*$. 
\end{itemize}
\end{theorem}

\bproof
It follows from \cite{LT2} that the analogue of (\ref{eq3.1}) holds for
all finite sets of vectors in $L$ and therefore the equivalence of (i)
and (ii) follows from Corollary \ref{cor2.5}.
\eproof

It is well known and easy to prove that if $X$ or $X^*$ has cotype 2
then $\GL$ and $\GL_2$ are equivalent for $X$. This leads to

\begin{theorem}
\label{thm3.3}
Let $X$ be a Banach space of cotype 2. The following statements are
equivalent:
\begin{itemize}
\item[(i)] $X$ has $\GL$.
\item[(ii)] $X$ has $S\GL(1,2)$.
\item[(iii)] $X$ has $S\GL(1,q)$ for some $q$, $1\le q<\infty$.
\item[(iv)]  There exists a Banach lattice $L$ with $X\subseteq L$ so
that $X$ is 2-concave in $L$.
\end{itemize}
\end{theorem}

\bproof
The equivalence between (i) and (ii) can be proved as \cite[Theorems 1.3
(iv) and 1.4 (i)]{CN}. Since $X$ is of cotype 2
$B(L_\infty,X)=\Pi_2(L_\infty,X)$, which together with an easy
application of Maurey's extension property \cite{M} shows that
$\Pi_q(Y,X)=\Pi_2(Y,X)$ for all $2\le q<\infty$ and all $Y$. This shows
the equivalence between (ii) and (iii). It follows directly from Theorem
\ref{thm2.7} that (ii) and (iv) are equivalent.
\eproof

As a corollary we obtain

\begin{corollary}
\label{cor3.4}
Any cotype 2 subspace of a Banach lattice of finite cotype has $\GL$.
\end{corollary}

Let us end this section with some results which relate the
$\GL$-property of a Banach space $X$ to compactness of absolutely
summing operators defined on $X$. Our first result is probably
well-known.

\begin{proposition}
\label{prop3.5}
Let $X$ and $Y$ be Banach spaces so that $X$ does not contain a subspace
isomorphic to $\ell_1$ and let $1\le p<\infty$. Then every $p$-summing
operator from $X$ to $Y$ is compact.
\end{proposition}

\bproof
Let $T\in\Pi_p(X,Y)$ and let $(x_n)\subseteq X$ with $\|x_n\|\le
1$. Since $X$ does not contain $\ell_1$ it follows from Rosenthal's
$\ell_1$-theorem \cite{LT1} that $(x_n)$ has a subsequence $(x_{n_k})$,
which is a weak Cauchy sequence. It now follows from a result of Pietsch
\cite{Pi} that $(Tx_{n_k})$ is norm convergent in $Y$.

Hence $T$ is compact.
\eproof

The next result can often be used to prove that a given concrete space
does not have $\GL$.

\begin{theorem}
\label{thm3.6}
Let $Y$ be a Banach space with $\GL$, and let $X$ be a quotient of $Y$.

If there exists an absolutely summing, non-compact operator from $X$ to
a Banach space $Z$ with the Radon-Nikodym property (RNP), then $Y$
contains a complemented subspace isomorphic to $\ell_1$.
\end{theorem}

\bproof
Let $Q\colon Y\to X$ be a quotient map of $Y$ onto $X$, let $Z$ be a
Banach space with the RNP and assume that there is a $T\in\Pi_1(X,Z)$
which is not compact. Hence there exists a sequence $(x_n)\subseteq X$
with $\|x_n\| <1$ for all $n\in\bN$ and an $\e>0$ so that
\begin{equation}
\label{eq3.3}
\| Tx_n-Tx_m\|\ge\e\quad\mbox{for all $n,m\in \bN$}.
\end{equation}
For every $n\in\bN$ we choose an $y_n\in Y$ with $\|y_n\|<1$ so that
$Qy_n=x_n$. Since $Y$ has $\GL$ there exist a measure $\mu$ and
operators $A\in B(Y,L_1(\mu))$, $B\in B(L_1(\mu),Z)$ so that $TQ=BA$. By
\cite{D} $B$ takes weak cauchy sequence into norm convergent ones and
therefore $(Ay_n)$ does not have any weak Cauchy subsequence and hence
by \cite{KP} there exist a subsequence $(Ay_{n_k})$ of $(Ay_n)$ which is
equivalent to the unit vector basis of $\ell_1$ and a bounded projection
$P$ of $L_1(\mu)$ onto $[Ay_{n_k}]$. Since $\ell_1$ has the lifting
property and $PA$ maps $Y$ onto $[Ay_{n_k}]$ it follows that if $U$ is
any isomorphism of $\ell_1$ onto $[Ay_{n_k}]$ then there exists a $V\in
B(\ell_1,Y)$ so that $U=PAV$. Clearly $V(\ell_1)$ is isomorphic to
$\ell_1$ and $VU^{-1}PA$ is a projection of $Y$ onto $V(\ell_1)$. 
\eproof

\begin{corollary}
\label{cor3.7}
Let $Y$ be a Banach space of finite cotype with $\GL$ and let
$X\subseteq Y$ be a subspace. Then every absolutely summing operator
from $X^*$ to a Banach space $Z$ with the RNP is compact.
\end{corollary}

\bproof
Since $Y$ is of finite cotype $\ell_1$ cannot be isomorphic to a
complemented subspace of $Y^*$ and hence the conclusion follows from
Theorem \ref{thm3.6}.
\eproof

\section{Some concluding remarks}
\label{sec4}

\setcounter{equation}{0}

The construction in Section \ref{sec1} gives rise to the hope that it
could be possible to develop a theory of lattice subspaces with the
so-called regular operators as morphisms, somewhat following the idea
from operator spaces. This turns out not to be possible, if one in
addition requires a reasonable duality theory. In this section we wish
to comment a little on these problems. We start with the following
definition:

\begin{definition}
\label{def4.1}
Let $L$ and $M$ be Banach lattices, $X\subseteq L$, $Y\subseteq M$ be
subspaces and \linebreak $T\in B(X,Y)$. $T$ is called $\ell_1$-regular
(respectively $\ell_\infty$-regular) if there is a constant $K\ge 1$ so
that for all finite sets $\{x_1,x_2,\dots,x_n\}\subseteq X$ we have
\begin{equation}
\label{eq4.1}
\|\sum^n_{j=1} |Tx_j| \| \le K\|\sum^n_{j=1} |x_j|\|
\end{equation}
(respectively
\begin{equation}
\label{eq4.2}
\|\bigvee^n_{j=1} |Tx_j| \| \le K\|\bigvee^n_{j=1} |x_j|\| \quad\big).
\end{equation}
\end{definition}

If $T\in B(L,M)$ then $\ell_\infty$-regularity of $T$ equals the usual
definition of a regular operator \cite{S}. It is easy to see that in
this case $T$ is regular if and only if it is $\ell_1$-regular. This
turns out not to be the case if $T$ is only defined on a subspace of a
Banach lattice, as the example below shows. Let us first state the
following lemma

\begin{lemma}
\label{lemma4.2}
Let $(\Delta,\cM,\mu)$ be a measure space and $X$ a subspace of
$L_1(\mu)$. If $L$ is a Banach lattice, then every $T\in B(X,L)$ is
$\ell_1$-regular. 
\end{lemma}

\bproof
If $T\in B(X,L)$ and $\{x_1,x_2,\dots,x_n\}\subseteq X$, then
\begin{equation}
\label{eq4.3}
\|\sum^n_{j=1} |Tx_j| \| \le \sum^n_{j=1} \|Tx_j\|\le \|T\| \sum^n_{j=1}
\|x_j\| = \|T\| \|\sum^n_{j=1} |x_j| \|.
\end{equation}
\eproof

\begin{example}
\label{ex4.3}
Let $(r_j)$ be the sequence of Rademacher functions in $L_1(0,1)$, put
$H=[r_j]$ and let $T\colon H\to \ell_2$ be the natural isomorphism.
By the lemma $T$ is $\ell_1$-regular, but $\|\bigvee^n_{j=1} |Tr_j| \|_2
= \|\sum^n_{j=1} Tr_j\|_2=\sqrt{n}$ and $\|\bigvee^n_{j=1} |r_j| \|_1=1$
which shows that $T$ is not $\ell_\infty$-regular.
\end{example}

If $X$ is a subspace of a Banach lattice $L$ then we can consider
$\ell_\infty\otimes X^*$ as a (non-closed) subspace of $(\ell_1\otimes_m
X)^*$ and define the norm $\alpha$ on $\ell_\infty\otimes X^*$ as the
restriction of the norm on $(\ell_1\otimes_m X)^*$. This $\alpha$ is
readily seen to be left invariant and hence Theorem \ref{thm1.7} gives
an embedding of $X^*$ into a Banach lattice $M$, so that $\alpha$
corresponds to the $m$-norm. With this construction one could try to
build up a theory of lattice subspaces using the $\ell_\infty$-regular
operators as morphisms. Unfortunately it will not lead to a reasonable
duality theory. Indeed, doing the above dualization twice we obtain an
embedding of $X^{**}$ into a Banach lattice, but the canonical embedding
of $X$ into $X^{**}$ need not be $\ell_\infty$-regular.

One of the main reasons for this obstacle is the difference between
$\ell_1$- and $\ell_\infty$-regularity.

In her Ph.D.-thesis L.B.\ McClaran \cite{MC} makes a thorough
investigation of subspaces and quotients of Banach lattices and has
succeeded in developing a theory for subspaces of quotients of Banach
lattices with the $\ell_\infty$- and $\ell_1$-regular operators as
morphisms.

There are rudiments of a duality theory in some of the results in the
previous sections and it is our belief that a theory of $\GL$-subspaces
of Banach lattices can be developed.

\vspace{1cm}

\noindent Department of Mathematics,\\ University of Missouri,\\ Columbia MO 65211,\\
pete@casazza.math.missouri.edu\\

\noindent Department of Mathematics and Computer Science,\\ Odense 
University,\\ Campusvej 55, DK-5230 Odense M, Denmark,\\
njn@imada.ou.dk

\end{document}